\documentclass{article}
\usepackage{amssymb}

\usepackage{amsmath}
\usepackage{graphicx}
\usepackage{color}

\linethickness{0.6pt}
\newsavebox{\fstring}
\sbox{\fstring}{{\color{blue}\circle*{3},$\rule{40pt}{0mm}$,\circle*{3},
$\rule{17pt}{0mm}$,\circle*{3},$\rule{16pt}{0mm}$,\circle*{3},
$\rule{16pt}{0mm}$,\circle*{3},$\rule{17pt}{0mm}$,\circle*{3},
$\rule{16pt}{0mm}$,\circle*{3},$\rule{17pt}{0mm}$,\circle*{3},
$\rule{16pt}{0mm}$,\circle*{3},$\rule{17pt}{0mm}$,\circle*{3}}}
\newsavebox{\fstringzm}
\sbox{\fstringzm}{$\rule{38pt}{0mm}$ -1
$\rule{7pt}{0mm}$ +2 $\rule{7pt}{0mm}$ -2
$\rule{7pt}{0mm}$ +2 $\rule{7pt}{0mm}$ -3
$\rule{7pt}{0mm}$ +4 $\rule{7pt}{0mm}$ -5
$\rule{7pt}{0mm}$ +6 $\rule{6pt}{0mm}$ -7}
\newsavebox{\fstringm}
\sbox{\fstringm}{-1 $\rule{25pt}{0mm}$ -1
$\rule{7pt}{0mm}$ +2 $\rule{7pt}{0mm}$ -2
$\rule{7pt}{0mm}$ +2 $\rule{7pt}{0mm}$ -3
$\rule{7pt}{0mm}$ +4 $\rule{7pt}{0mm}$ -5
$\rule{7pt}{0mm}$ +6 $\rule{6pt}{0mm}$ -7}
\newsavebox{\fsstring}
\sbox{\fsstring}{{\color{blue}\circle*{3},$\rule{40pt}{0mm}$,\circle*{3},
$\rule{17pt}{0mm}$,\circle*{3},$\rule{16pt}{0mm}$,\circle*{3},
$\rule{16pt}{0mm}$,\circle*{3},$\rule{17pt}{0mm}$,\circle*{3},
$\rule{16pt}{0mm}$,\circle*{3}}}
\newsavebox{\fsstringm}
\sbox{\fsstringm}{-1 $\rule{25pt}{0mm}$ -1
$\rule{7pt}{0mm}$ +2 $\rule{7pt}{0mm}$ -2
$\rule{7pt}{0mm}$ +2 $\rule{7pt}{0mm}$ -3
$\rule{7pt}{0mm}$ +4}
\newsavebox{\fstrings}
\sbox{\fstrings}{{\color{blue}\circle*{3},
$\rule{16pt}{0mm}$,\circle*{3},$\rule{17pt}{0mm}$,\circle*{3},
$\rule{17pt}{0mm}$,\circle*{3},$\rule{16pt}{0mm}$,\circle*{3},
$\rule{16pt}{0mm}$,\circle*{3},$\rule{17pt}{0mm}$,\circle*{3},
$\rule{16pt}{0mm}$,\circle*{3},$\rule{17pt}{0mm}$,\circle*{3},
$\rule{16pt}{0mm}$,\circle*{3},$\rule{17pt}{0mm}$,\circle*{3}}}
\newsavebox{\fstringsm}
\sbox{\fstringsm}{+1
$\rule{6pt}{0mm}$ -1 $\rule{6pt}{0mm}$ +1
$\rule{6pt}{0mm}$ -1 $\rule{6pt}{0mm}$ +2
$\rule{6pt}{0mm}$ -3 $\rule{6pt}{0mm}$ +3
$\rule{6pt}{0mm}$ -4 $\rule{6pt}{0mm}$ +6
$\rule{6pt}{0mm}$ -7 $\rule{6pt}{0mm}$ +8}
\newsavebox{\fsstrings}
\sbox{\fsstrings}{{\color{blue}\circle*{3},
$\rule{16pt}{0mm}$,\circle*{3},$\rule{17pt}{0mm}$,\circle*{3},
$\rule{17pt}{0mm}$,\circle*{3},$\rule{16pt}{0mm}$,\circle*{3},
$\rule{16pt}{0mm}$,\circle*{3},$\rule{17pt}{0mm}$,\circle*{3},
$\rule{16pt}{0mm}$,\circle*{3},$\rule{17pt}{0mm}$,\circle*{3},
$\rule{16pt}{0mm}$,\circle*{3}}}
\newsavebox{\fsstringsm}
\sbox{\fsstringsm}{+1
$\rule{6pt}{0mm}$ -1 $\rule{6pt}{0mm}$ +1
$\rule{6pt}{0mm}$ -1 $\rule{6pt}{0mm}$ +2
$\rule{6pt}{0mm}$ -3 $\rule{6pt}{0mm}$ +3
$\rule{6pt}{0mm}$ -4 $\rule{6pt}{0mm}$ +6
$\rule{6pt}{0mm}$ -7}
\newsavebox{\fssstrings}
\sbox{\fssstrings}{{\color{blue}\circle*{3},
$\rule{16pt}{0mm}$,\circle*{3},$\rule{17pt}{0mm}$,\circle*{3},
$\rule{17pt}{0mm}$,\circle*{3},$\rule{16pt}{0mm}$,\circle*{3},
$\rule{16pt}{0mm}$,\circle*{3},$\rule{17pt}{0mm}$,\circle*{3},
$\rule{16pt}{0mm}$,\circle*{3},$\rule{17pt}{0mm}$,\circle*{3}}}
\newsavebox{\fssstringsm}
\sbox{\fssstringsm}{+1
$\rule{6pt}{0mm}$ -1 $\rule{6pt}{0mm}$ +1
$\rule{6pt}{0mm}$ -1 $\rule{6pt}{0mm}$ +2
$\rule{6pt}{0mm}$ -3 $\rule{6pt}{0mm}$ +3
$\rule{6pt}{0mm}$ -4 $\rule{6pt}{0mm}$ +6}
\newsavebox{\fsssstrings}
\sbox{\fsssstrings}{{\color{blue}\circle*{3},
$\rule{16pt}{0mm}$,\circle*{3},$\rule{17pt}{0mm}$,\circle*{3},
$\rule{17pt}{0mm}$,\circle*{3},$\rule{16pt}{0mm}$,\circle*{3},
$\rule{16pt}{0mm}$,\circle*{3}}}
\newsavebox{\fsssstringsm}
\sbox{\fsssstringsm}{+1
$\rule{6pt}{0mm}$ -1 $\rule{6pt}{0mm}$ +1
$\rule{6pt}{0mm}$ -1 $\rule{6pt}{0mm}$ +2
$\rule{6pt}{0mm}$ -3}
\newsavebox{\fssssstrings}
\sbox{\fssssstrings}{{\color{blue}\circle*{3},
$\rule{16pt}{0mm}$,\circle*{3},$\rule{17pt}{0mm}$,\circle*{3},
$\rule{17pt}{0mm}$,\circle*{3}}}
\newsavebox{\fssssstringsm}
\sbox{\fssssstringsm}{+1
$\rule{6pt}{0mm}$ -1 $\rule{6pt}{0mm}$ +1
$\rule{6pt}{0mm}$ -1}
\newsavebox{\fbranch}
\sbox{\fbranch}{{\Huge $\ast$}
$\rule{25pt}{0mm}$ {\Huge $\ast$} $\rule{26pt}{0mm}$ {\Huge $\ast$}
$\rule{26pt}{0mm}$ {\Huge $\ast$} $\rule{73pt}{0mm}$ {\Huge $\ast$}}
\newsavebox{\fbranchm}
\sbox{\fbranchm}{-4
$\rule{30pt}{0mm}$ -3 $\rule{28pt}{0mm}$  -2
$\rule{31pt}{0mm}$  -1 $\rule{82pt}{0mm}$  -1}
\newsavebox{\fbranchp}
\sbox{\fbranchp}{+4
$\rule{26pt}{0mm}$ +3 $\rule{26pt}{0mm}$  +2
$\rule{26pt}{0mm}$  +1 $\rule{75pt}{0mm}$  +1}
\newsavebox{\fbrranch}
\sbox{\fbrranch}{{\Huge $\ast$} $\rule{26pt}{0mm}$ {\Huge $\ast$}
$\rule{26pt}{0mm}$ {\Huge $\ast$} $\rule{73pt}{0mm}$ {\Huge $\ast$}}
\newsavebox{\fbrranchm}
\sbox{\fbrranchm}{-3 $\rule{28pt}{0mm}$  -2
$\rule{31pt}{0mm}$  -1 $\rule{82pt}{0mm}$  -1}
\newsavebox{\fbrranchp}
\sbox{\fbrranchp}{ +3 $\rule{26pt}{0mm}$  +2
$\rule{26pt}{0mm}$  +1 $\rule{75pt}{0mm}$  +1}
\newsavebox{\fbrrranch}
\sbox{\fbrrranch}{{\Huge $\ast$} $\rule{26pt}{0mm}$ {\Huge $\ast$}
$\rule{26pt}{0mm}$ {\Huge $\ast$} $\rule{73pt}{0mm}$ {\Huge $\ast$}}
\newsavebox{\fbrrranchm}
\sbox{\fbrrranchm}{-3 $\rule{28pt}{0mm}$  -2
$\rule{31pt}{0mm}$  -1 $\rule{82pt}{0mm}$  -1}
\newsavebox{\fbrrranchp}
\sbox{\fbrrranchp}{ +3 $\rule{26pt}{0mm}$  +2
$\rule{26pt}{0mm}$  +1 $\rule{75pt}{0mm}$  +1}
\newsavebox{\fbranchs}
\sbox{\fbranchs}{$\rule{13pt}{0mm}$ {\Huge $\ast$} $\rule{26pt}{0mm}$ {\Huge $\ast$}
$\rule{26pt}{0mm}$ {\Huge $\ast$} $\rule{26pt}{0mm}$ {\Huge $\ast$}
$\rule{26pt}{0mm}$ {\Huge $\ast$}}
\newsavebox{\fbranchms}
\sbox{\fbranchms}{$\rule{13pt}{0mm}$ -5 $\rule{28pt}{0mm}$  -4
$\rule{31pt}{0mm}$  -2 $\rule{31pt}{0mm}$  -2 $\rule{31pt}{0mm}$  -1}
\newsavebox{\fbranchps}
\sbox{\fbranchps}{$\rule{13pt}{0mm}$ +5 $\rule{26pt}{0mm}$  +4
$\rule{26pt}{0mm}$  +2 $\rule{27pt}{0mm}$  +2 $\rule{27pt}{0mm}$  +1}
\newsavebox{\fbrranchs}
\sbox{\fbrranchs}{$\rule{13pt}{0mm}$ {\Huge $\ast$} $\rule{26pt}{0mm}$ {\Huge $\ast$}
$\rule{26pt}{0mm}$ {\Huge $\ast$}}
\newsavebox{\fbrranchms}
\sbox{\fbrranchms}{$\rule{13pt}{0mm}$ -2 $\rule{31pt}{0mm}$  -2 $\rule{31pt}{0mm}$  -1}
\newsavebox{\fbrranchps}
\sbox{\fbrranchps}{$\rule{13pt}{0mm}$ +2 $\rule{27pt}{0mm}$  +2 $\rule{27pt}{0mm}$  +1}
\newsavebox{\fbrrranchs}
\sbox{\fbrrranchs}{{\Huge $\ast$} $\rule{26pt}{0mm}$ {\Huge $\ast$}}
\newsavebox{\fbrrranchms}
\sbox{\fbrrranchms}{ -2 $\rule{31pt}{0mm}$  -1}
\newsavebox{\fbrrranchps}
\sbox{\fbrrranchps}{ +2 $\rule{27pt}{0mm}$  +1}
\newtheorem{theorem}{Theorem}

\newtheorem{corollary}[theorem]{Corollary}

\newtheorem{example}[theorem]{Example}

\newtheorem{proposition}[theorem]{Proposition}

\input{tcilatex}

\begin{document}

\title{\textbf{{\Large {On a property of branching coefficients }\\
for affine Lie algebras}}}
\author{Mikhail Ilyin \\
Theoretical Department, SPb State University,\\
198904, Sankt-Petersburg, Russia \\
[2ex] Petr Kulish\thanks{
Supported by RFFI grant N 06-01-00451 }\\
Sankt-Petersburg Department of\\ Steklov Institute of Mathematics\\
Fontanka 27, 191023, Sankt-Petersburg, Russia \\
[5mm] Vladimir Lyakhovsky \thanks{
Supported by RFFI grant N 06-01-00451 and the National Project
RNP.2.1.1.1112 }\\
Theoretical Department, SPb State University,\\
198904, Sankt-Petersburg, Russia \\
e-mail:lyakh1507@nm.ru \\
}
\maketitle

\begin{abstract}
It is demonstrated that decompositions of integrable highest weight modules
of a simple Lie algebra with respect to its reductive subalgebra obey the
set of algebraic relations leading to the recursive properties for the
corresponding branching coefficients. These properties are encoded in the
special element $\Gamma _{\frak{g} \supset \frak{a}}$ of the formal algebra $%
\mathcal{E}_{\frak{a}}$ that describes the injection and is called the fan.
In the simplest case, when $\frak{a} = \frak{h} \left( \frak{g} \right)$,
the recursion procedure generates the weight diagram of a module $L_{\frak{g}%
}$. When applied to a reduction of highest weight modules the recursion
described by the fan provides a highly effective tool to obtain the explicit
values of branching coefficients.
\end{abstract}

\section{ Introduction}

We consider integrable modules $L^{\mu }$ of affine Lie algebra $\frak{g}$
with the highest weight $\mu $ and the reduced modules $L_{\frak{g\downarrow
a}}^{\mu }$ with respect to a reductive subalgebra $\frak{a}\subset \frak{g}$%
. In particular when the Cartan subalgebra $\frak{a}= \frak{h}$ is studied,
the branching coefficients indicate the dimensions of the weight subspaces
and thus describe the module. String functions and branching coefficients of
the affine Lie algebra pairs (e.g. $A_n^{(1)} \subset A_{n-p-1}^{(1)} \oplus
A_p^{(1)}$) arise in the computation of the local state probabilities for
solvable models on square lattice \cite{DJKMO}. Irreducible highest weight
modules with dominant integral weights appear also in application of the
quantum inverse scattering method \cite{LD} where solvable spin chains are
studied in the framework of the AdS/CFT correspondence conjecture of the
super-string theory (see \cite{KAZA,BE} and references therein).

There are different ways to find branching coefficients. One can use the BGG
resolution \cite{BGG} (for Kac-Moody algebras the algorithm is described in
\cite{Kac,Wak1}), the Schure function series \cite{FauKing}, the BRST
cohomology \cite{Hwang}, Kac-Peterson formulas \cite{Kac} or the
combinatorial methods applied in \cite{FeigJimbo}. We want to obtain the
recursive formulas for weight multiplicities and branching coefficients
using the purely algebraic approach. From the Weyl-Kac character formula
\cite{Kac}
\begin{equation}
\mathrm{ch}L^{\mu }\left( \frak{g}\right) =\frac{\sum_{w\in W}\epsilon
(w)e^{w\circ (\mu +\rho )-\rho }}{\prod_{\alpha \in \Delta
^{+}}(1-e^{-\alpha })^{\mathrm{mult}(\alpha )}},
\end{equation}
we derive the special set of relations for
branching coefficients. These relations can be used both to construct a
representation and to reduce it with respect to a subalgebra $\frak{a}%
\subset \frak{g}$ that is to find the corresponding branching rules. Each
relation of the set deals with a finite collection of weights. Among them it
is always possible to fix the lowest one (with respect to the natural
ordering of weights induced by basic roots). Thus it is possible to use the
relations of the set as recurrent relations for branching coefficients. It
is demonstrated that branching is governed by a certain system of weights
(called ''the fan of injection'') that depends only on the algebra and
the injection morphism and can be used to decompose the highest weight
modules.

For finite dimensional classical Lie algebras the case of regular injections
was considered in \cite{LF} where the recurrent relations were constructed
using the properties of the Kostant-Heckman partition function. The same
method was used for regular injections of affine Lie algebras \cite{LDu}. In
this study we present a different approach and find that for any reductive
subalgebra $\frak{a}$ of an affine Lie algebra $\frak{g}$ such that $\frak{h}%
_{\frak{a}}^{\ast }\subset \frak{h}_{\frak{g}}^{\ast }$ and $\frak{h}_{%
\overset{\circ }{\frak{a}}}^{\ast }\subset \frak{h}_{\overset{\circ }{\frak{g%
}}}^{\ast }$\ the branching coefficients obey the set of properties that
give rise to recurrent relations. The latter provide a compact and
effective method to construct the corresponding branching rules. The results
are illustrated by examples.

\section{Basic definitions and relations.}

Consider the affine Lie algebras $\frak{g}$ and $\frak{a}$ with the
underlying finite-dimensional subalgebras $\overset{\circ }{\frak{g}}$ and $%
\overset{\circ }{\frak{a}}$ and an injection $\frak{a}\longrightarrow \frak{g%
}$ such that $\frak{a}$ is a reductive subalgebra $\frak{a\subset g}$ with
correlated root spaces: $\frak{h}_{\frak{a}}^{\ast }\subset \frak{h}_{\frak{g%
}}^{\ast }$ and $\frak{h}_{\overset{\circ }{\frak{a}}}^{\ast }\subset \frak{h%
}_{\overset{\circ }{\frak{g}}}^{\ast }$\ .

The following notation will also be used:

$L^{\mu }$\ $\left( L_{\frak{a}}^{\nu }\right) $\ -- the integrable module
of $\frak{g}$ with the highest weight $\mu $\ ; (resp. integrable $\frak{a}$%
-module with the highest weight $\nu $ );

$r$ , $\left( r_{\frak{a}}\right) $ -- the rank of the algebra $\frak{g}$ $%
\left( \text{resp. }\frak{a}\right) $ ;

$\Delta $ $\left( \Delta _{\frak{a}}\right) $-- the root system; $\Delta
^{+} $ $\left( \text{resp. }\Delta _{\frak{a}}^{+}\right) $-- the positive
root system (of $\frak{g}$ and $\frak{a}$ respectively);

$\mathrm{mult}\left( \alpha \right) $ $\left( \mathrm{mult}_{\frak{a}}\left(
\alpha \right) \right) $ -- the multiplicity of the root $\alpha$ in $\Delta
$ (resp. in $\left( \Delta _{\frak{a}}\right) $);

$\overset{\circ }{\Delta }$ , $\left( \overset{\circ }{\Delta _{\frak{a}}}%
\right) $ -- the finite root system of the subalgebra $\overset{\circ }{%
\frak{g}}$ (resp. $\overset{\circ }{\frak{a}}$);

$\mathcal{N}^{\mu }$ , $\left( \mathcal{N}_{\frak{a}}^{\nu }\right) $ -- the
weight diagram of $L^{\mu }$ $\left( \text{resp. }L_{\frak{a}}^{\nu }\right)
$ ;

$W$ , $\left( W_{\frak{a}}\right) $-- the corresponding Weyl group;

$C$ , $\left( C_{\frak{a}}\right) $-- the fundamental Weyl chamber;

$\rho $\ , $\left( \rho _{\frak{a}}\right) $\ -- the Weyl vector;

$\epsilon \left( w\right) :=\det \left( w\right) $ ;

$\alpha _{i}$ , $\left( \alpha _{\left( \frak{a}\right) j}\right) $ -- the $i
$-th (resp. $j$-th) basic root for $\frak{g}$ $\left( \text{resp. }\frak{a}%
\right) $; $i=0,\ldots ,r$ ,\ \ $\left( j=0,\ldots ,r_{\frak{a}}\right) $;

$\delta $ -- the imaginary root of $\frak{g}$ (and of $\frak{a}$ if any);

$\alpha _{i}^{\vee }$ , $\left( \alpha _{\left( \frak{a}\right) j}^{\vee
}\right) $-- the basic coroot for $\frak{g}$ $\left( \text{resp. }\frak{a}%
\right) $ , $i=0,\ldots ,r$ ;\ \ $\left( j=0,\ldots ,r_{\frak{a}}\right) $;

$\overset{\circ }{\xi }$ , $\overset{\circ }{\xi _{\left( \frak{a}\right) }}$
-- the finite (classical) part of the weight $\xi \in P$ , $\left( \text{%
resp. }\xi _{\left( \frak{a}\right) }\in P_{\frak{a}}\right) $\ ;

$\lambda =\left( \overset{\circ }{\lambda };k;n\right) $ -- the
decomposition of an affine weight indicating the finite part $\overset{\circ
}{\lambda }$, level $k$ and grade $n$\ .

$P$ $\left( \text{resp. } P_{\frak{a}}\right) $ \ -- the weight lattice;

$M \left( \text{resp. }M_{\frak{a}}\right) :=$

\noindent $=\left\{
\begin{array}{c}
\sum_{i=1}^{r}\mathbf{Z}\alpha _{i}^{\vee }\text{ }\left( \text{resp. }%
\sum_{i=1}^{r}\mathbf{Z}\alpha _{\left( \frak{a}\right) i}^{\vee }\right)
\text{for untwisted algebras or }A_{2r}^{\left( 2\right) }, \\
\sum_{i=1}^{r}\mathbf{Z}\alpha _{i}\text{ }\left( \text{resp. }\sum_{i=1}^{r}%
\mathbf{Z}\alpha _{\left( \frak{a}\right) i}\right) \text{for }A_{r}^{\left(
u\geq 2\right) }\text{ and }A\neq A_{2r}^{\left( 2\right) },
\end{array}
\right\} ;$

$\mathcal{E}$\ , $\left( \mathcal{E}_{\frak{a}}\right) $-- the group algebra
of the group $P$ (resp. $P_{\frak{a}} $);

$\Theta _{\lambda }:=e^{-\frac{\left| \lambda \right| ^{2}}{2k}\delta
}\sum\limits_{\alpha \in M}e^{t_{\alpha }\circ \lambda }$ -- the classical
theta-function;

$\Theta _{\left( \frak{a}\right) \nu }:=e^{-\frac{\left| \nu \right| ^{2}}{%
2k_{\frak{a}}}\delta }\sum\limits_{\beta \in M_{\frak{a}}}e^{t_{\beta }\circ
\nu }$;\newline
notice that when the injection is considered the level $k_{\frak{a}}$ must
be correlated with the corresponding rescaling of roots;

$A_{\lambda }:=\sum\limits_{s\in \overset{\circ }{W}}\epsilon (s)\Theta
_{s\circ \lambda }$ $\left( \text{resp. }A_{\left( \frak{a}\right) \nu
}:=\sum\limits_{s\in \overset{\circ }{W_{\frak{a}}}}\epsilon (s)\Theta
_{\left( \frak{a}\right) s\circ \nu }\right) $;

$\Psi ^{\left( \mu \right) }:=e^{\frac{\left| \mu +\rho \right| ^{2}}{2k}%
\delta \ -\ \rho }A_{\mu +\rho }=e^{\frac{\left| \mu +\rho \right| ^{2}}{2k}%
\delta \ -\ \rho }\sum\limits_{s\in \overset{\circ }{W}}\epsilon (s)\Theta
_{s\circ \left( \mu +\rho \right) }=$

\noindent $=\sum\limits_{w\in W}\epsilon (w)e^{w\circ (\mu +\rho )-\rho }$
-- the singular weight element for the $\frak{g}$-module $L^{\mu }$;

$\Psi _{\left( \frak{a}\right) }^{\left( \nu \right) }:=e^{\frac{\left| \nu
+\rho _{_{\frak{a}}}\right| ^{2}}{2k_{\frak{a}}}\delta \ -\ \rho _{_{\frak{a}%
}}}A_{\left( \frak{a}\right) \nu +\rho _{_{\frak{a}}}}=e^{\frac{\left| \nu
+\rho _{_{\frak{a}}}\right| ^{2}}{2k_{\frak{a}}}\delta \ -\ \rho _{_{\frak{a}%
}}}\sum\limits_{s\in \overset{\circ }{W_{\frak{a}}}}\epsilon (s)\Theta
_{\left( \frak{a}\right) s\circ \left( \nu +\rho _{_{\frak{a}}}\right) }=$

\noindent $=\sum\limits_{w\in W_{\frak{a}}}\epsilon (w)e^{w\circ (\nu +\rho
_{_{\frak{a}}})-\rho _{_{\frak{a}}}}$ -- the corresponding singular weight
element for the $\frak{a}$-module $L_{\frak{a}}^{\nu }$;

$\widehat{\Psi ^{\left( \mu \right) }}$ $\left( \widehat{\Psi _{\left( \frak{%
a}\right) }^{\left( \nu \right) }}\right) $ -- the set of singular weights $%
\xi \in P$ $\left( \text{resp. }\in P_{\frak{a}}\right) $ for the module $%
L^{\mu }$ $\left( \text{resp. }L_{\frak{a}}^{\nu }\right) $ with the
coordinates $\left( \overset{\circ }{\xi },k,n,\epsilon \left( w\left( \xi
\right) \right) \right) \mid _{\xi =w\left( \xi \right) \circ (\mu +\rho
)-\rho },$ (resp. $\left( \overset{\circ }{\xi },k,n,\epsilon \left(
w_{a}\left( \xi \right) \right) \right) \mid _{\xi =w_{a}\left( \xi \right)
\circ (\nu +\rho _{a})-\rho _{a}}$ ), (this set is similar to $P_{\mathrm{%
nice}}^{\prime }\left( \mu \right) $ in \cite{Wak1})

$m_{\xi }^{\left( \mu \right) }$ , $\left( m_{\xi }^{\left( \nu \right)
}\right) $ -- the multiplicity of the weight $\xi \in P$ \ $\left( \text{%
resp. }\in P_{\frak{a}}\right) $ in the module $L^{\mu }$ , (resp. $\xi \in
L_{\frak{a}}^{\nu } $);

$ch\left( L^{\mu }\right) $ $\left( \text{resp. }ch\left( L_{\frak{a}}^{\nu
}\right) \right) $-- the formal character of $L^{\mu }$ $\left( \text{resp. }%
L_{\frak{a}}^{\nu }\right) $;

$ch\left( L^{\mu }\right) =\frac{\sum_{w\in W}\epsilon (w)e^{w\circ (\mu
+\rho )-\rho }}{\prod_{\alpha \in \Delta ^{+}}\left( 1-e^{-\alpha }\right) ^{%
\mathrm{{mult}\left( \alpha \right) }}}=\frac{\Psi ^{\left( \mu \right) }}{%
\Psi ^{\left( 0\right) }}$ -- the Weyl-Kac formula.

$R:=\prod_{\alpha \in \Delta ^{+}}\left( 1-e^{-\alpha }\right) ^{\mathrm{{%
mult}\left( \alpha \right) }}=\Psi ^{\left( 0\right) }\quad $

\noindent $\left( \text{resp. }R_{\frak{a}}:=\prod_{\alpha \in \Delta _{%
\frak{a}}^{+}}\left( 1-e^{-\alpha }\right) ^{\mathrm{mult}_{\frak{a}}\mathrm{%
\left( \alpha \right) }}=\Psi _{\frak{a}}^{\left( 0\right) }\right) $-- the
denominator.

\section{Anomalous multiplicities and recurrent relations}

For the injection $\frak{a\longrightarrow g}$ consider the reduced module
\begin{equation}
L_{\frak{g}\downarrow \frak{a}}^{\mu }=\bigoplus\limits_{\nu \in P_{\frak{a}%
}^{+}}b_{\nu }^{\left( \mu \right) }L_{\frak{a}}^{\nu }  \label{reduction}
\end{equation}
with the branching coefficients $b_{\nu }^{\left( \mu \right) }$. The
character reduction
\begin{equation}
\pi _{\frak{a}}\circ \left( \mathrm{ch}L_{\frak{g}}^{\mu }\right) =\sum_{\nu
\in P_{\frak{a}}^{+}}b_{\nu }^{\left( \mu \right) }\mathrm{ch}L_{\frak{a}%
}^{\nu }  \label{char-dec}
\end{equation}
involves the projection operator $\pi _{\frak{a}}:P\longrightarrow P_{\frak{a%
}}$.

The denominator identity can be applied to redress the relation (\ref
{char-dec}),
\begin{equation}
\frac{\pi _{\frak{a}}\circ \left( \sum_{w\in W}\epsilon (w)e^{w\circ (\mu
+\rho )-\rho }\right) }{\pi _{\frak{a}}\circ \left( \prod_{\alpha \in \Delta
^{+}}\left( 1-e^{-\alpha }\right) ^{\mathrm{{mult}\left( \alpha \right) }%
}\right) }=\sum_{\nu \in P_{\frak{a}}^{+}}b_{\nu }^{\left( \mu \right) }%
\frac{\sum\limits_{w\in W_{\frak{a}}}\epsilon (w)e^{w\circ (\nu +\rho _{_{%
\frak{a}}})-\rho _{_{\frak{a}}}}}{\prod_{\beta \in \Delta _{\frak{a}%
}^{+}}\left( 1-e^{-\beta }\right) ^{\mathrm{mult}_{\frak{a}}\mathrm{\left(
\beta \right) }}},  \label{intermediate}
\end{equation}
and to rewrite it as
\begin{equation}
\begin{array}{c}
\pi _{\frak{a}}\circ \left( \sum_{w\in W}\epsilon (w)e^{w\circ (\mu +\rho
)-\rho }\right) = \\
\rule{0mm}{10mm}=\frac{\pi _{\frak{a}}\circ \left( \prod_{\alpha \in \Delta
^{+}}\left( 1-e^{-\alpha }\right) ^{\mathrm{{mult}\left( \alpha \right) }%
}\right) }{\prod_{\beta \in \Delta _{\frak{a}}^{+}}\left( 1-e^{-\beta
}\right) ^{\mathrm{mult}_{\frak{a}}\mathrm{\left( \beta \right) }}}\sum_{\nu
\in P_{\frak{a}}^{+}}b_{\nu }^{\left( \mu \right) }\sum\limits_{w\in W_{%
\frak{a}}}\epsilon (w)e^{w\circ (\nu +\rho _{_{\frak{a}}})-\rho _{_{\frak{a}%
}}},
\end{array}
\label{origin}
\end{equation}
For the trivial $\frak{g}$-module $L^{0}$ with $\mu =0$ we have
\begin{equation*}
\frac{\pi _{\frak{a}}\circ \left( \sum_{w\in W}\epsilon (w)e^{w\circ \rho
-\rho }\right) }{\pi _{\frak{a}}\circ \left( \prod_{\alpha \in \Delta
^{+}}\left( 1-e^{-\alpha }\right) ^{\mathrm{{mult}\left( \alpha \right) }%
}\right) }=\frac{\sum\limits_{w\in W_{\frak{a}}}\epsilon (w)e^{w\circ \rho
_{_{\frak{a}}}-\rho _{_{\frak{a}}}}}{\prod_{\beta \in \Delta _{\frak{a}%
}^{+}}\left( 1-e^{-\beta }\right) ^{\mathrm{mult}_{\frak{a}}\mathrm{\left(
\beta \right) }}},
\end{equation*}
\begin{equation*}
\frac{\pi _{\frak{a}}\circ \left( \prod_{\alpha \in \Delta ^{+}}\left(
1-e^{-\alpha }\right) ^{\mathrm{{mult}\left( \alpha \right) }}\right) }{%
\prod_{\beta \in \Delta _{\frak{a}}^{+}}\left( 1-e^{-\beta }\right) ^{%
\mathrm{mult}_{\frak{a}}\mathrm{\left( \beta \right) }}}=\frac{\pi _{\frak{a}%
}\circ \left( \sum_{w\in W}\epsilon (w)e^{w\circ \rho -\rho }\right) }{%
\sum\limits_{w\in W_{\frak{a}}}\epsilon (w)e^{w\circ \rho _{_{\frak{a}%
}}-\rho _{_{\frak{a}}}}}.
\end{equation*}
The relation (\ref{origin}) takes the form
\begin{equation}
\begin{array}{l}
\pi _{\frak{a}}\circ \left( \sum_{w\in W}\epsilon (w)e^{w\circ (\mu +\rho
)-\rho }\right) \sum\limits_{w\in W_{\frak{a}}}\epsilon (w)e^{w\circ \rho
_{_{\frak{a}}}-\rho _{_{\frak{a}}}}= \\
\rule{0mm}{8mm}=\pi _{\frak{a}}\circ \left( \sum_{w\in W}\epsilon
(w)e^{w\circ \rho -\rho }\right) \sum_{\nu \in P_{\frak{a}}^{+}}b_{\nu
}^{\left( \mu \right) }\sum\limits_{w\in W_{\frak{a}}}\epsilon (w)e^{w\circ
(\nu +\rho _{_{\frak{a}}})-\rho _{_{\frak{a}}}}.
\label{origin-reformulated}
\end{array}
\end{equation}
Consider the expression $\sum_{\nu \in P_{\frak{a}}^{+}}b_{\nu
}^{\left( \mu \right) }\Psi _{\left( \frak{%
a}\right) }^{\left( \nu \right) }$ and introduce the
numbers $k_{\lambda }^{\left( \mu \right) }$, called \emph{the anomalous
branching coefficients}, -- the multiplicities of submodules $L_{\frak{a}%
}^{\nu }$ times the determinants $\epsilon \left( w\right) $ contained in $%
\Psi _{\left( \frak{a}\right) }^{\left( \nu \right) }$.
\begin{equation}
\sum_{\nu \in P_{\frak{a}}^+}b_{\nu }^{\left( \mu \right) }\Psi _{\left( \frak{%
a}\right) }^{\left( \nu \right) }=\sum_{\lambda \in P_{\frak{a}}}k_{\lambda
}^{\left( \mu \right) }e^{\lambda }  \label{anom-br-coeff}
\end{equation}
In these terms the expression (\ref{origin-reformulated}) reads
\begin{equation}
\begin{array}{c}
\pi _{\frak{a}}\circ \left( \sum_{w\in W}\epsilon (w)e^{w\circ (\mu +\rho
)-\rho }\right) \sum\limits_{w\in W_{\frak{a}}}\epsilon (w)e^{w\circ \rho
_{_{\frak{a}}}-\rho _{_{\frak{a}}}}= \\
\rule{0mm}{8mm}=\pi _{\frak{a}}\circ \left( \sum_{w\in W}\epsilon
(w)e^{w\circ \rho -\rho }\right) \sum_{\xi \in P_{\frak{a}}}k_{\xi }^{\left(
\mu \right) }e^{\xi },
\end{array}
\end{equation}
or
\begin{equation}
\begin{array}{c}
\sum_{w\in W,v\in W_{\frak{a}}}\epsilon (w)\epsilon (v)e^{\pi _{\frak{a}%
}\circ \left( w\circ (\mu +\rho )-\rho \right) +v\circ \rho _{_{\frak{a}%
}}-\rho _{_{\frak{a}}}}= \\
\rule{0mm}{7mm}=\sum_{\xi \in P_{\frak{a}}}\sum_{w\in W}\epsilon (w)e^{\pi _{%
\frak{a}}\circ \left( w\circ \rho -\rho \right) +\xi }k_{\xi }^{\left( \mu
\right) }.
\end{array}
\end{equation}
Thus we have proved the statement:

\begin{proposition}
Let $L^{\mu }$ be the integrable highest weight module of $\frak{g}$, $\frak{%
a}\subset \frak{g}$, $\frak{h}_{\frak{a}}\subset \frak{h}_{\frak{g}}$, $%
\frak{h}_{\frak{a}}^{\ast }\subset \frak{h}_{\frak{g}}^{\ast }$ and $\pi _{%
\frak{a}}$ -- a projection $P\longrightarrow P_{\frak{a}}$ -- then for any
point $\xi \in P_{\frak{a}}$ the following relation holds:
\begin{equation}
\sum_{w\in W}\epsilon (w)k_{\xi -\pi _{\frak{a}}\circ \left( w\circ \rho
-\rho \right) }^{\left( \mu \right) }=\sum_{w\in W,v\in W_{\frak{a}%
}}\epsilon (w)\epsilon (v)\delta _{\pi _{\frak{a}}\circ \left( w\circ (\mu
+\rho )-\rho \right) ,\xi +\rho _{_{\frak{a}}}-v\circ \rho _{_{\frak{a}}}}.
\label{property-with-star}
\end{equation}
\end{proposition}

This can be rearranged to produce a recurrent relation for the anomalous
multiplicities,
\begin{equation}
k_{\xi }^{\left( \mu \right) }=-\sum_{w\in W\setminus e}\epsilon (w)k_{\xi
-\pi _{\frak{a}}\circ \left( w\circ \rho -\rho \right) }^{\left( \mu \right)
}+\sum_{w\in W,v\in W_{\frak{a}}}\epsilon (w)\epsilon (v)\delta _{\pi _{%
\frak{a}}\circ \left( w\circ (\mu +\rho )-\rho \right) ,\xi +\rho _{_{\frak{a%
}}}-v\circ \rho _{_{\frak{a}}}}.  \label{branching-with-star}
\end{equation}
This formula can be applied to find the branching coefficients $b_{\nu
}^{\left( \mu \right) }$ due to the fact that being restricted to the
fundamental Weyl chamber $\left( C_{\frak{a}}\right) $ the anomalous
branching coefficients coincide with the branching coefficients
\begin{equation*}
k_{\xi }^{\left( \mu \right) }=b_{\xi }^{\left( \mu \right) }\qquad \mathrm{%
for\quad }\xi \in C_{\frak{a}}.
\end{equation*}

The relation (\ref{branching-with-star}) contains the standard system of
shifts, $\xi \longrightarrow \xi -\pi _{\frak{a}}\circ \left( w\circ \rho
-\rho \right) $, defined by the singular weights of the trivial module, the
corresponding element of the algebra $\mathcal{E}$ being $\Psi ^{\left(
0\right) }=e^{\frac{\left| \rho \right| ^{2}}{2k}\delta \ -\ \rho }$ $\times
\sum\limits_{s\in \overset{\circ }{W}}\epsilon (s)\Theta _{s\circ \left(
\rho \right) }$\ . At the same time the second term in the r.h.s. contains
the summation in both $W$ and $W_{\frak{a}}$\ . Below we demonstrate that
this relation can be simplified by introducing the different system of
shifts.

Let us return to the relation (\ref{origin}). The conditions $\frak{%
a\longrightarrow g}$ and $\frak{h}_{\frak{a}}\subset \frak{h}_{\frak{g}}$
guarantee the inclusion $\Delta _{\frak{a}}^{+}\subset \Delta ^{+}$. Thus
the first factor in the r. h. s. being an element of $\mathcal{E}$ can be
written as
\begin{eqnarray*}
\frac{\pi _{\frak{a}}\circ \left( \prod_{\alpha \in \Delta ^{+}}\left(
1-e^{-\alpha }\right) ^{\mathrm{{mult}\left( \alpha \right) }}\right) }{%
\prod_{\beta \in \Delta _{\frak{a}}^{+}}\left( 1-e^{-\beta }\right) ^{%
\mathrm{mult}_{\frak{a}}\mathrm{\left( \beta \right) }}} &=&\prod_{\alpha
\in \left( \pi _{\frak{a}}\circ \Delta ^{+}\right) }\left( 1-e^{-\alpha
}\right) ^{\mathrm{{mult}\left( \alpha \right) -{mult}}_{\frak{a}}\mathrm{%
\left( \alpha \right) }}= \\
&=&-\sum_{\gamma \in P_{\frak{a}}}s\left( \gamma \right) e^{-\gamma }.
\end{eqnarray*}
For the coefficient function $s\left( \gamma \right) $ define $\Phi _{\frak{a%
}\subset \frak{g}}\subset P_{\frak{a}}$ as its carrier:
\begin{equation}
\Phi _{\frak{a}\subset \frak{g}}=\left\{ \gamma \in P_{\frak{a}}\mid s\left(
\gamma \right) \neq 0\right\} ;  \label{phi-d}
\end{equation}
\begin{equation}
\prod_{\alpha \in \left( \pi _{\frak{a}}\circ \Delta ^{+}\right) }\left(
1-e^{-\alpha }\right) ^{\mathrm{{mult}\left( \alpha \right) -{mult}}_{\frak{a%
}}\mathrm{\left( \alpha \right) }}=-\sum_{\gamma \in \Phi _{\frak{a}\subset
\frak{g}}}s\left( \gamma \right) e^{-\gamma }.  \label{fan-d}
\end{equation}
When the second factor in the r.h.s. of (\ref{origin}) is also decomposed we
obtain the relation
\begin{eqnarray*}
\sum_{w\in W}\epsilon (w)e^{\pi _{\frak{a}}\circ \left( w\circ (\mu +\rho
)-\rho \right) } &=&-\sum_{\gamma \in \Phi _{\frak{a}\subset \frak{g}%
}}s\left( \gamma \right) e^{-\gamma }\sum_{\lambda \in P_{\frak{a}%
}}k_{\lambda }^{\left( \mu \right) }e^{\lambda } \\
&=&-\sum_{\gamma \in \Phi _{\frak{a}\subset \frak{g}}}\sum_{\lambda \in P_{%
\frak{a}}}s\left( \gamma \right) k_{\lambda }^{\left( \mu \right)
}e^{\lambda -\gamma }
\end{eqnarray*}
and the new property
\begin{equation}
\sum_{w\in W}\epsilon \left( w\right) \delta _{\xi ,\pi _{\frak{a}}\circ
\left( w\circ (\mu +\rho )-\rho \right) }+\sum_{\gamma \in \Phi _{\frak{a}%
\subset \frak{g}}}s\left( \gamma \right) k_{\xi +\gamma }^{\left( \mu
\right) }=0;\qquad \xi \in P_{\frak{a}}.  \label{property}
\end{equation}

Thus the following statement is true:

\begin{proposition}
Let $L^{\mu }$ be the integrable highest weight module of $\frak{g}$, $\frak{%
a}\subset \frak{g}$, $\frak{h}_{\frak{a}}\subset \frak{h}_{\frak{g}}$, $%
\frak{h}_{\frak{a}}^{\ast }\subset \frak{h}_{\frak{g}}^{\ast }$ and $\pi _{%
\frak{a}}$ -- a projection $P\longrightarrow P_{\frak{a}}$ then for any
vector $\xi \in P_{\frak{a}}$ the sum $\sum_{\gamma \in \Phi _{\frak{a}%
\subset \frak{g}}}-s\left( \gamma \right) k_{\xi +\gamma }^{\left( \mu
\right) }$ is equal to the anomalous multiplicity of the weight $\xi $ in
the module $\pi _{\frak{a}}\circ L_{\frak{g}}^{\mu }$.
\end{proposition}

This property (\ref{property}) also produces recurrent relations for the
anomalous multiplicities. Returning to the relation (\ref{fan-d}) we see
that $\Phi _{\frak{a}\subset \frak{g}}$ contains vectors with nonnegative
grade and is a subset in the carrier of the singular weights element $\Psi
^{\left( \mu \right) }=e^{\frac{\left| \mu +\rho \right| ^{2}}{2k}\delta \
-\ \rho }\sum\limits_{s\in \overset{\circ }{W}}\epsilon (s)\Theta _{s\circ
\left( \mu +\rho \right) }$. In each grade the set $\Phi _{\frak{a}\subset
\frak{g}}$ has finite number of vectors \cite{Kac}. In particular
\begin{equation*}
\#\left( \Phi _{\frak{a}\subset \frak{g}}\right) _{n=0}=\#\overset{\circ }{W}%
.
\end{equation*}
\ In $\left( \Phi _{\frak{a}\subset \frak{g}}\right) _{n=0}$ let $\gamma
_{0} $ be the lowest vector with respect to the natural ordering in $%
\overset{\circ }{\Delta _{\frak{a}}}$ . Decomposing the defining relation (%
\ref{fan-d}),
\begin{equation}
\prod_{\alpha \in \left( \pi _{\frak{a}}\circ \Delta ^{+}\right) }\left(
1-e^{-\alpha }\right) ^{\mathrm{{mult}\left( \alpha \right) -{mult}}_{\frak{a%
}}\mathrm{\left( \alpha \right) }}=-s\left( \gamma _{0}\right) e^{-\gamma
_{0}}-\sum_{\gamma \in \Phi _{\frak{a}\subset \frak{g}}\setminus \gamma
_{0}}s\left( \gamma \right) e^{-\gamma },  \label{pre-fan-def}
\end{equation}
in (\ref{property}) we obtain

\begin{equation}
k_{\xi }^{\left( \mu \right) }=-\frac{1}{s\left( \gamma _{0}\right) }\left(
\sum_{w\in W}\epsilon \left( w\right) \delta _{\xi ,\pi _{\frak{a}}\circ
\left( w\circ (\mu +\rho )-\rho \right) +\gamma _{0}}+\sum_{\gamma \in
\Gamma _{\frak{a}\subset \frak{g}}}s\left( \gamma +\gamma _{0}\right) k_{\xi
+\gamma }^{\left( \mu \right) }\right)   \label{recurrent relation}
\end{equation}
where the set
\begin{equation}
\Gamma _{\frak{a}\subset \frak{g}}=\left\{ \xi -\gamma _{0}|\xi \in \Phi _{%
\frak{a}\subset \frak{g}}\right\} \setminus \left\{ 0\right\} .
\label{fan-defined}
\end{equation}
was introduced called \emph{the fan of the injection} $\frak{a}\subset \frak{%
g}$. The equality (\ref{recurrent relation}) can be considered as a
recurrent relation for anomalous branching coefficients $k_{\xi }^{\left(
\mu \right) }$. Contrary to the relation (\ref{branching-with-star}) here
only the summation over $W$ is applied.

When $r=r_{\frak{a}}$ the positive roots $\Delta _{\frak{a}}^{+}$ can always
be chosen so that $\gamma _{0}=0$, the relation (\ref{pre-fan-def})
indicates that $s\left( \gamma _{0}\right) =-1$. Thus in this special case
\begin{equation*}
\Gamma _{\frak{a}\subset \frak{g}}=\Phi _{\frak{a}\subset \frak{g}}\setminus
\left\{ 0\right\} ,
\end{equation*}
and the recurrent relation acquires the form
\begin{equation}
k_{\xi }^{\left( \mu \right) }=\sum_{\gamma \in \Gamma _{\frak{a}\subset
\frak{g}}}s\left( \gamma \right) k_{\xi +\gamma }^{\left( \mu \right)
}+\sum_{w\in W}\epsilon \left( w\right) \delta _{\xi ,\pi _{\frak{a}}\circ
\left( w\circ \left( \mu +\rho \right) -\rho \right) }.
\label{recursion-rel}
\end{equation}

Comments:

\begin{enumerate}
\item  The sets $\Phi _{\frak{a}\subset \frak{g}}$ and $\Gamma _{\frak{a}%
\subset \frak{g}}$ do not depend on the representation $L^{\mu }\left( \frak{%
g}\right) $ and describe the injection of the subalgebra $\frak{a}$ into the
algebra $\frak{g}$.

\item  Let the set of singular weights of the projected module $\pi _{\frak{a%
}}\circ L^{\mu }$ be constructed. Then the sets $\Psi ^{\left( 0\right) }$
and $\Psi _{\left( \frak{a}\right) }^{\left( 0\right) }$ define the
anomalous branching coefficients for the reduced module $L_{\frak{%
g\downarrow a}}^{\mu }$ by means of relation (\ref{branching-with-star}).
The same information can be obtained using the fan $\Gamma _{\frak{a}\subset
\frak{g}}$ via the relations (\ref{recurrent relation}) or (\ref
{recursion-rel}).

\item  The set of branching coefficients $\left\{ b_{\nu }^{\left( \mu
\right) }\right\} $ is the subset of the anomalous branching coefficients $%
\left\{ k_{\xi }^{\left( \mu \right) }\right\} $ :
\begin{equation*}
\left\{ b_{\nu }^{\left( \mu \right) }\mid \nu \in P_{\frak{a}}^{+}\right\}
=\left\{ k_{\xi }^{\left( \mu \right) }\mid \xi \in \overline{C_{\frak{a}}}%
\right\} .
\end{equation*}
Thus the recurrent relations (\ref{branching-with-star}),(\ref{recurrent
relation}) and (\ref{recursion-rel}) supply us with the branching
coefficients as well.
\end{enumerate}

Let us apply the obtained results to the case where $\frak{a}$ is a Cartan
subalgebra of $\frak{g}$ , $\frak{a}=\frak{h}_{\frak{g}}$. Then the Weyl
group $W_{\frak{a}}$ and the projector $\pi _{\frak{a}}$ are trivial and in
the formulas (\ref{branching-with-star}) and (\ref{recursion-rel}) the
anomalous coefficient $k_{\xi }^{\left( \mu \right) }$ is the multiplicity
of the only singular weight that contributes to the element $\Psi _{\frak{h}%
_{\frak{g}}}^{\left( \xi \right) }$ -- the highest weight of the $\frak{h}_{%
\frak{g}}$-submodule. This means that here the number $k_{\xi }^{\left( \mu
\right) }$ is always nonnegative and coincides with the multiplicity $m_{\xi
}^{\left( \mu \right) }$ of the weight $\xi $ in the module $L^{\mu }\left(
\frak{g}\right) $. The relations (\ref{branching-with-star}) and (\ref
{recursion-rel}) directly lead to

\begin{corollary}
In the integrable highest weight module $L^{\mu }\left( \frak{g}\right) $ of
an affine Lie algebra $\frak{g}$ the multiplicity $m_{\xi }^{\left( \mu
\right) }$ of the weight $\xi $ (considered as a numerical function on $P_{%
\frak{g}}$) obeys the relation
\begin{equation}
m_{\xi }^{\left( \mu \right) }=-\sum_{w\in W\setminus e}\epsilon (w)m_{\xi
-\left( w\circ \rho -\rho \right) }^{\left( \mu \right) }+\sum_{w\in
W}\epsilon (w)\delta _{\left( w\circ (\mu +\rho )-\rho \right) ,\xi }.
\label{weights-recursion}
\end{equation}
\end{corollary}
In implicit form this relation can be found for affine Lie algebras in
\cite{Kac}
(Ch.11, the second formula in Ex. 11.14) and for finite dimensional
algebras in
\cite{Ful}.
Usually the truncated formula (without the second term in the r. h. s.)
is presented as the recurrent relation for the multiplicities of weights
(see for example \cite{Burbaki}, Ch. VIII, Sect. 9.3). From our
point of view it is highly important to deal with the recurrent relation
in its full form.
The reason is that the relation (\ref{weights-recursion}) can be applied
for any reducible module and is valid in any domain of $P$, not only
inside the diagram
$\mathcal{N}^{\mu } \setminus \mu$ with the single highest weight $\mu$.

This relation gives the possibility to construct recursively the module $%
L^{\mu }\left( \frak{g}\right) $ provided the elements
$\Psi ^{\left( \mu \right) }$ and $\Psi ^{\left( 0\right) }$ are known.

\newpage

\section{ Examples}

\setcounter{theorem}{0}

\begin{example}
Consider the finite dimensional Lie algebras $A_{2}\subset g_{2}$. The root
system $\Delta $ is generated by the simple roots $\alpha _{1}$ (the long)
and $\alpha _{2}$ (the short) with the angle $\frac{5\pi }{6}$ between them.
\begin{equation}
\Delta ^{+}=\left\{ \alpha _{1},\alpha _{2},\alpha _{1}+\alpha _{2},\alpha
_{1}+2\alpha _{2},\alpha _{1}+3\alpha _{2},2\alpha _{1}+3\alpha _{2}\right\}
\label{pos-roots-g2}
\end{equation}
\begin{equation}
\Delta _{sl(3)}^{+}=\left\{ \alpha _{1},\alpha _{1}+3\alpha _{2},2\alpha
_{1}+3\alpha _{2}\right\}   \label{pos-roots-sl3}
\end{equation}
From (\ref{phi-d}), (\ref{fan-d}) and (\ref{fan-defined}) we obtain
\begin{equation}
\Phi _{A_{2}\subset g_{2}}=\left\{ 0,\alpha _{2},\alpha _{1}+\alpha
_{2},\alpha _{1}+3\alpha _{2},2\alpha _{1}+3\alpha _{2},2\alpha _{1}+4\alpha
_{2}\right\}   \label{pre-fan-sl3-g2}
\end{equation}
\begin{equation}
\Gamma _{A_{2}\subset g_{2}}=\left\{ \alpha _{2},\alpha _{1}+\alpha
_{2},\alpha _{1}+3\alpha _{2},2\alpha _{1}+3\alpha _{2},2\alpha _{1}+4\alpha
_{2}\right\}   \label{fan-sl3-g2}
\end{equation}
Consider the adjoint module $L^{2\alpha _{1}+3\alpha _{2}}$. Its singular
weights are
\begin{eqnarray}
&&\left\{ 2\alpha _{1}+3\alpha _{2},3\alpha _{2},-\alpha _{1}+2\alpha
_{2},-4\alpha _{2},-\alpha _{1}-6\alpha _{2},-8\alpha _{1}-12\alpha
_{2},-8\alpha _{1}-13\alpha _{2},\right.   \notag \\
&&\left. -6\alpha _{1}-13\alpha _{2},-5\alpha _{1}-12\alpha _{2},-6\alpha
_{1}-6\alpha _{2},-5\alpha _{1}-4\alpha _{2},2\alpha _{1}+2\alpha
_{2}\right\}   \label{anom-adj-g2}
\end{eqnarray}
to each of them corresponds the Weyl transformation $w\left( \psi \right) $
and the value $\epsilon \left( w\right) $:
\begin{equation}
\left\{ \epsilon \left( w\left( \psi \right) \right) \right\} =\left\{
+1,-1,+1,+1,-1,-1,+1,-1,+1,-1,+1,-1\right\}   \label{sfunc-sl3-g2}
\end{equation}
In the closure of the fundamental chamber $\overline{C_{\frak{a}}}$ the
relation (\ref{recursion-rel}) defines three nonzero branching coefficients
\begin{equation}
b_{2\alpha _{1}+3\alpha _{2}}^{\left( \mu \right) }=+1,\quad b_{\alpha
_{1}+2\alpha _{2}}^{\left( \mu \right) }=+1,\quad b_{\alpha _{1}+\alpha
_{2}}^{\left( \mu \right) }=+1.  \label{branchcoeff-sl3-g2}
\end{equation}
corresponding to the adjoint and two fundamental submodules of $sl(3)$ in
the decomposition $L_{\downarrow sl(3)}^{2\alpha _{1}+3\alpha _{2}}$ (
Notice that we need the singular weight $2\alpha _{1}+2\alpha _{2}$ with $%
s\left( 2\alpha _{1}+2\alpha _{2}\right) =-1$ to be used in the above
calculations.)
\end{example}

\begin{example}
Consider the special injection of the algebra $B_{1}$ into $A_{2}$. Let $%
\alpha _{1}$ and $\alpha _{2}$ be the simple roots of $A_{2}$,
\begin{equation}
\Delta ^{+}=\left\{ \alpha _{1},\alpha _{2},\alpha _{1}+\alpha _{2}\right\}
\label{canonical-pos-roots-sl3}
\end{equation}
The only positive root of $B_{1}$ is
\begin{equation}
\Delta _{B_{1}}^{+}=\left\{ \beta :=\frac{1}{2}\alpha _{1}\right\}
\label{pos-roots-so3}
\end{equation}
From (\ref{phi-d}) it follows that
\begin{equation}
\Phi _{B_{1}\subset A_{2}}=\left\{ 0,\alpha _{1},-\frac{1}{2}\alpha _{1},%
\frac{1}{2}\alpha _{1}\right\} =\left\{ 0,2\beta ,-\beta ,\beta \right\}
\label{pre-fan-so3-sl3}
\end{equation}
For these vectors the function $s_{B_{1}\subset A_{2}}$ has the values
\begin{equation}
s_{B_{1}\subset A_{2}}=\left\{ -1,+1,+1,-1\right\}
\label{sign-func-so3-sl3}
\end{equation}
The minimal vector $\gamma _{0}$
\begin{eqnarray}
\gamma _{0}^{B_{1}\subset A_{2}} &=&-\beta   \notag \\
s_{B_{1}\subset A_{2}}\left( -\beta \right)  &=&+1.
\label{min-vect-so3-sl3}
\end{eqnarray}
The fan is formed by eliminating $\gamma _{0}$ from $\Phi _{B_{1}\subset
A_{2}}$ and shifting the remaining vectors by $-\gamma _{0}$:
\begin{eqnarray}
\Gamma _{B_{1}\subset A_{2}} &=&\left\{ \beta ,2\beta ,3\beta \right\}
\notag \\
s_{B_{1}\subset A_{2}}\left( \gamma +\gamma _{0}\right)  &=&\left\{
-,\,-,\,+\,\right\} ;\qquad \gamma \in \Gamma _{B_{1}\subset A_{2}}
\label{fan-so3-sl3}
\end{eqnarray}
Consider the module $L^{\alpha _{1}+\alpha _{2}}$. Its singular
weights are
\begin{equation}
\left\{ \alpha _{1}+\alpha _{2},-\alpha _{1}+\alpha _{2},-3\alpha
_{1}-\alpha _{2},-3\alpha _{1},-\alpha _{1}-3\alpha _{2},\alpha _{1}-\alpha
_{2}\right\}   \label{anom-adj-sl3}
\end{equation}
with the values
\begin{equation}
\left\{ \epsilon \left( \xi \right) \right\} =\left\{ +1,\qquad -1,\qquad
+1,\qquad -1,\qquad +1,\qquad -1\right\}   \label{sign-func-sl3-again}
\end{equation}
We rewrite their projections to $P_{B_{1}}$ in terms of $\beta $:
\begin{eqnarray}
&&\left\{ -6\beta ,-4\beta ,-4\beta ,0,0,+2\beta \right\}   \notag \\
\left\{ \epsilon \left( \xi \right) \right\}  &=&\left\{
+1,\,-1,\,+1,\,-1,\,+1,\,-1\,\right\}   \label{anom-adj-sl3-shift}
\end{eqnarray}
In the closure of the fundamental chamber $\left( \overline{C_{\frak{a}}}%
\right) $ the relation (\ref{recurrent relation}) defines two nonzero
branching coefficients
\begin{equation}
b_{2\beta }^{\left( \mu \right) }=+1,\quad b_{\beta }^{\left( \mu \right)
}=+1.  \label{branchcoeff-so3-sl3}
\end{equation}
corresponding to the submodule of the adjoint subrepresentation
and the 5-dimensional spin 2
submodule of $B_{1}$ in the reduced module $L_{\downarrow B_{1}}^{\alpha
_{1}+\alpha _{2}}$ (Notice that the singular vector ''$0$'' with $s\left(
0\right) =-1$ has the multiplicity 2.)
\end{example}

\begin{example}
For the affine algebra $A_{2}^{(1)}$ consider the twisted subalgebra $%
A_{2}^{(2)}$. For the level $k$ sublattice $P_{k}$ introduce the normalized
basic vectors $\left\{ e_{1},e_{2},e_{3}\right\} $ with $\mid e_{j}\mid
_{j=1,2,3}=1$\ and $\delta $ \ with $\mid \delta \mid =0$ .\ For $A_{2}^{(1)}
$ we fix the simple roots
\begin{equation*}
\alpha _{1}=e_{1}-e_{2};\quad \alpha _{2}=e_{2}-e_{3};\quad \alpha
_{0}=\delta -e_{1}+e_{3};
\end{equation*}
The positive roots are as follows:
\begin{equation*}
\Delta ^{+}=\left\{
\begin{array}{c}
\alpha _{j}+l\delta ;\quad j=1,2,3;\quad l\in \mathbf{Z}_{\geq 0} \\
-\alpha _{j}+p\delta ;\quad j=1,2,3;\quad p\in \mathbf{Z}_{>0} \\
p\delta ;\quad mult\left( p\delta \right) =2;\quad p\in \mathbf{Z}_{>0}
\end{array}
\right\} ,
\end{equation*}
the classical positive roots being
\begin{equation*}
\overset{\circ }{\Delta }^{+}=\left\{ \alpha _{1}=e_{1}-e_{2};\quad \alpha
_{2}=e_{2}-e_{3};\quad \alpha _{3}=e_{1}-e_{3};\right\} .
\end{equation*}
The fundamental weights
\begin{equation*}
\omega _{1}=\frac{1}{3}\left( 2e_{1}-e_{2}-e_{3}\right) +k;\quad \omega _{2}=%
\frac{1}{3}\left( e_{2}+e_{1}-2e_{3}\right) +k;\quad \omega _{0}=k;
\end{equation*}
and the Weyl vector
\begin{equation*}
\rho =(\alpha _{1}+\alpha _{2},3,0).
\end{equation*}
The Weyl group is generated by the classical reflections
\begin{equation*}
s_{\alpha _{1}},s_{\alpha _{2}}
\end{equation*}
and (in accord with $M=\sum_{i=1}^{r}\mathbf{Z}\alpha _{i}^{\vee }$ for
untwisted algebras) the translations
\begin{equation*}
t_{\alpha _{1}},t_{\alpha _{2}}.
\end{equation*}

Consider the module $L^{\omega _{0}}$. Notice that to obtain the branching
rules we need only the projected singular element $\pi _{\frak{a}}\circ \Psi
^{\left( \omega _{0}\right) }$ of this module and the set $\Gamma
_{A_{2}^{(2)}\subset A_{2}^{(1)}}$ and do not need any other properties of
the module itself. Let us describe the element $\Psi ^{\left( \omega
_{0}\right) }$ by the set $\widehat{\Psi ^{\left( \omega _{0}\right) }}$ of
singular weights of the module $L^{\omega _{0}}$:
\begin{equation*}
\left\{ \left( \lambda _{1},\lambda _{2},\lambda _{3},n,\epsilon \left(
w\right) \right) |\lambda _{i}\in \mathbf{Z},n\in \mathbf{Z}_{\leq
0},\epsilon \left( w\right) =\pm 1\right\} ,
\end{equation*}
(the level is always $k=1$ and is not indicated). Then for $n>-10$ the set $%
\widehat{\Psi ^{\left( \omega _{0}\right) }}$ contains the following 54
vectors:
\begin{eqnarray*}
&&\widehat{\Psi ^{\left( \omega _{0}\right) }}=\left\{ {}\right.
(0,0,0,0,-1),(-1,1,0,0,1),(0,-1,1,0,1),(-1,-1,2,0,-1), \\
&&(-2,1,1,0,-1),(-2,0,2,0,1),(2,-3,1,-2,-1),(-4,0,4,-2,-1), \\
&&(-4,3,1,-2,1),(-1,-3,4,-2,1),(-1,3,-2,-2,-1),(2,0,-2,-2,1), \\
&&(-5,1,4,-3,1),(-4,-1,5,-3,1),(0,-4,4,-3,-1),(-2,-3,5,-3,-1), \\
&&(-4,4,0,-3,-1),(-5,3,2,-3,-1),(3,-3,0,-3,1),(2,-4,2,-3,1), \\
&&(0,3,-3,-3,1),(-2,4,-2,-3,1),(3,-1,-2,-3,-1),(2,1,-3,-3,-1), \\
&&(-5,-1,6,-5,-1),(-6,1,5,-5,-1),(0,-5,5,-5,1),(-2,-4,6,-5,1), \\
&&(-5,5,0,-5,1),(-6,4,2,-5,1),(4,-4,0,-5,-1),(3,-5,2,-5,-1), \\
&&(0,4,-4,-5,-1),(-2,5,-3,-5,-1),(3,1,-4,-5,1),(4,-1,-3,-5,1), \\
&&(-6,0,6,-6,1),(-1,-5,6,-6,-1),(-6,5,1,-6,-1),(4,-5,1,-6,1), \\
&&(-1,5,-4,-6,1),(4,0,-4,-6,-1),(-4,-4,8,-9,-1),(-5,-3,8,-9,1), \\
&&(-8,4,4,-9,-1),(-8,3,5,-9,1),(3,-7,4,-9,1),(2,-7,5,-9,-1), \\
&&(-4,7,-3,-9,1),(-5,7,-2,-9,-1),(6,-3,-3,-9,-1),(6,-4,-2,-9,1), \\
&&(3,3,-6,-9,-1),(2,4,-6,-9,1),\ldots \left. {}\right\} .
\end{eqnarray*}
\newpage The $\pi _{\frak{a}}$ projection leads to the following set of
vectors:

\noindent
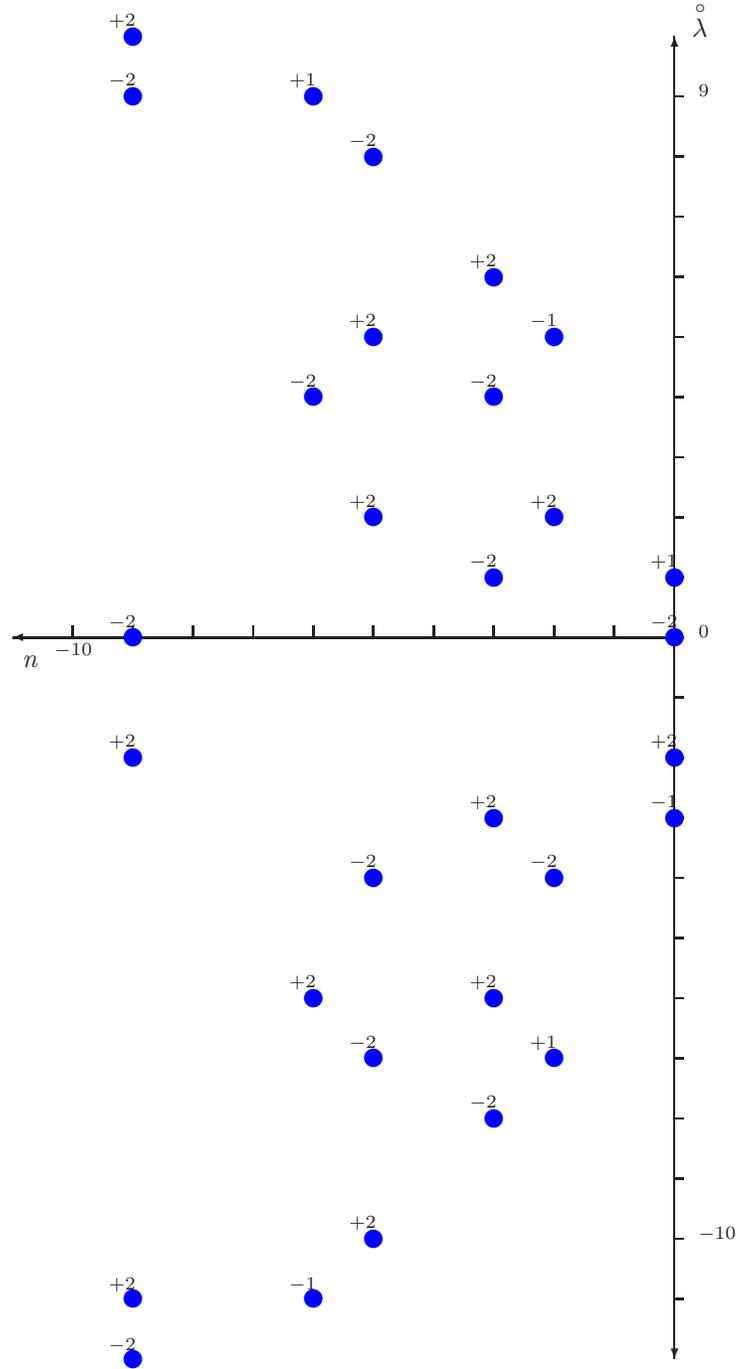
\begin{picture}(120.00,220.00)(-18.00,-13.00)
\put(120.00,100.00){\vector(0,50){100.00}}
\put(120.00,100.00){\vector(0,-50){120.00}}
\put(120.00,100.00){\vector(-1,0){110.00}}
\multiput(20.00,100.00)(10.00,00.00){10}{\line(0,1){2.00}}
\put(17.00,97.00){\scriptsize $-10$}
\put(10.00,95.00){ n}
\put(124.00,100.00){\scriptsize $0$}
\multiput(120.00,-10.00)(00.00,10.00){21}{\line(1,0){1.50}}
\put(124.00,190.00){\scriptsize $9$}
\put(123.00,200.00){$\overset{\circ }{\lambda }$}
\put(124.00,00.00){\scriptsize $-10$}
\put(120.00,100.00){{\color{blue}\circle*{3}}}
\put(116.00,101.80){\scriptsize $-2$}
\put(120.00,110.00){{\color{blue}\circle*{3}}}
\put(116.00,111.80){\scriptsize $+1$}
\put(120.00,80.00){{\color{blue}\circle*{3}}}
\put(116.00,81.80){\scriptsize $+2$}
\put(120.00,70.00){{\color{blue}\circle*{3}}}
\put(116.00,71.80){\scriptsize $-1$}
\put(100.00,120.00){{\color{blue}\circle*{3}}}
\put(96.00,121.80){\scriptsize $+2$}
\put(100.00,150.00){{\color{blue}\circle*{3}}}
\put(96.00,151.80){\scriptsize $-1$}
\put(100.00,60.00){{\color{blue}\circle*{3}}}
\put(96.00,61.80){\scriptsize $-2$}
\put(100.00,30.00){{\color{blue}\circle*{3}}}
\put(96.00,31.80){\scriptsize $+1$}
\put(90.00,160.00){{\color{blue}\circle*{3}}}
\put(86.00,161.80){\scriptsize $+2$}
\put(90.00,140.00){{\color{blue}\circle*{3}}}
\put(86.00,141.80){\scriptsize $-2$}
\put(90.00,110.00){{\color{blue}\circle*{3}}}
\put(86.00,111.80){\scriptsize $-2$}
\put(90.00,70.00){{\color{blue}\circle*{3}}}
\put(86.00,71.80){\scriptsize $+2$}
\put(90.00,40.00){{\color{blue}\circle*{3}}}
\put(86.00,41.80){\scriptsize $+2$}
\put(90.00,20.00){{\color{blue}\circle*{3}}}
\put(86.00,21.80){\scriptsize $-2$}
\put(70.00,180.00){{\color{blue}\circle*{3}}}
\put(66.00,181.80){\scriptsize $-2$}
\put(70.00,150.00){{\color{blue}\circle*{3}}}
\put(66.00,151.80){\scriptsize $+2$}
\put(70.00,120.00){{\color{blue}\circle*{3}}}
\put(66.00,121.80){\scriptsize $+2$}
\put(70.00,60.00){{\color{blue}\circle*{3}}}
\put(66.00,61.80){\scriptsize $-2$}
\put(70.00,30.00){{\color{blue}\circle*{3}}}
\put(66.00,31.80){\scriptsize $-2$}
\put(70.00,00.00){{\color{blue}\circle*{3}}}
\put(66.00,01.80){\scriptsize $+2$}
\put(60.00,190.00){{\color{blue}\circle*{3}}}
\put(56.00,191.80){\scriptsize $+1$}
\put(60.00,140.00){{\color{blue}\circle*{3}}}
\put(56.00,141.80){\scriptsize $-2$}
\put(60.00,40.00){{\color{blue}\circle*{3}}}
\put(56.00,41.80){\scriptsize $+2$}
\put(60.00,-10.00){{\color{blue}\circle*{3}}}
\put(56.00,-08.50){\scriptsize $-1$}
\put(30.00,200.00){{\color{blue}\circle*{3}}}
\put(26.00,201.80){\scriptsize $+2$}
\put(30.00,190.00){{\color{blue}\circle*{3}}}
\put(26.00,191.80){\scriptsize $-2$}
\put(30.00,100.00){{\color{blue}\circle*{3}}}
\put(26.00,101.80){\scriptsize $-2$}
\put(30.00,80.00){{\color{blue}\circle*{3}}}
\put(26.00,81.80){\scriptsize $+2$}
\put(30.00,-10.00){{\color{blue}\circle*{3}}}
\put(26.00,-08.50){\scriptsize $+2$}
\put(30.00,-20.00){{\color{blue}\circle*{3}}}
\put(26.00,-18.50){\scriptsize $-2$}
\end{picture}
\newpage \noindent In the figure the grade values increase to the right. The
vertical line unit is the basic root vector $\beta $ of $A_{2}^{(2)}$.

For the subalgebra $A_{2}^{(2)}$ the basic roots are
\begin{equation*}
\beta =\left( 1,0,0\right) ;\quad \beta _{0}=\delta -2\beta =\left(
-2,0,1\right) ;
\end{equation*}
with
\begin{equation*}
\theta =2\beta
\end{equation*}
and the normalization
\begin{equation}
\left| \beta \right| ^{2}=1,\quad \left| \beta _{0}\right| ^{2}=4,\quad
\left( \beta _{0},\beta \right) =-2.
\end{equation}
The fundamental weights
\begin{equation*}
\omega _{1}=1/2\beta +k=\left( 1/2,1,0\right) ,\quad \omega _{0}=2k=\left(
0,2,0\right)
\end{equation*}
and the Weyl vector
\begin{equation*}
\rho =1/2\beta +3k=(1/2,3,0).
\end{equation*}
The positive roots are
\begin{equation*}
\Delta _{A_{2}^{(2)}}^{+}=\left\{
\begin{array}{c}
\beta +n\delta ,\pm 2\beta +\left( 2n+1\right) \delta ;\qquad n\in \mathbf{Z}%
_{\geq 0} \\
-\beta +m\delta ;\quad m\delta \qquad m\in \mathbf{Z}_{>0}
\end{array}
\right\}
\end{equation*}
and have the multiplicity one. The Weyl group $W_{A_{2}^{(2)}}$ is generated
by the classical reflection $s_{\beta }$ and the translations $t\in
T_{A_{2}^{(2)}}\subset W_{A_{2}^{(2)}}$ along the coroot $\alpha _{0}^{\vee
}=1/2\delta -\beta =\left( -1,0,1/2\right) $: $T_{A_{2}^{(2)}}=\left\{
t_{l\alpha _{0}^{\vee }},l\in \mathbf{Z}\right\} $.

The injection $A_{2}^{(2)}\longrightarrow A_{2}^{(1)}$ is governed by its
classical part -- the special injection $B_{1}\longrightarrow A_{2}$. The
latter means that when we construct the subset $\Delta _{A_{2}^{(2)}}$ in
the root space of $A_{2}^{(1)}$ the roots in $\Delta _{A_{2}^{(2)}}$ are
scaled:
\begin{equation}
\beta =\alpha _{1}/2,\qquad K_{A_{2}^{(2)}}=2K_{A_{2}^{(1)}}.
\label{scale-3}
\end{equation}
(So in the modules $L^{\mu }$ of the level $k$ the $A_{2}^{(2)}$-submodules
have the level $2k$.)

According to (\ref{fan-d}) the set $\Phi _{A_{2}^{(2)}\subset A_{2}^{(1)}}$
is defined by the opposite vectors in the nonzero components of the element $%
\prod_{\alpha \in \left( \pi _{A_{2}^{(2)}}\circ \Delta ^{+}\right) }\left(
1-e^{-\alpha }\right) ^{\mathrm{{mult}\left( \alpha \right) -{mult}}_{\frak{a%
}}\mathrm{\left( \alpha \right) }}$. Taking into account the scaling (\ref
{scale-3}) we obtain the element
\begin{equation*}
\left( 1-e^{-\beta }\right) \left( 1-e^{2\beta }\right) \times
\prod\limits_{\varkappa =\pm 1}\prod\limits_{n=1}^{\infty }\left(
1-e^{\varkappa 2\beta +2n\delta }\right) \left( 1-e^{\varkappa \beta
+n\delta }\right) \prod\limits_{m=1}^{\infty }\left( 1-e^{+m\delta }\right)
\end{equation*}
generated by the set $\Phi _{A_{2}^{(2)}\subset A_{2}^{(1)}}$. Here the
lowest vector $\gamma _{0}$ is $-\beta $ and the fan is
\begin{equation*}
\Gamma _{A_{2}^{(2)}\subset A_{2}^{(1)}}=\left\{ \xi +\beta |\xi \in \Phi
_{A_{2}^{(2)}\subset A_{2}^{(1)}}\right\} \setminus \left\{ 0\right\} .
\end{equation*}
The structure of the fan can be illustrated by the following figure
presenting the vectors $\gamma \in \Gamma _{A_{2}^{(2)}\subset A_{2}^{(1)}}$
with the grade $n\leq 10$ and their multiplicities $s\left( \gamma \right) $.

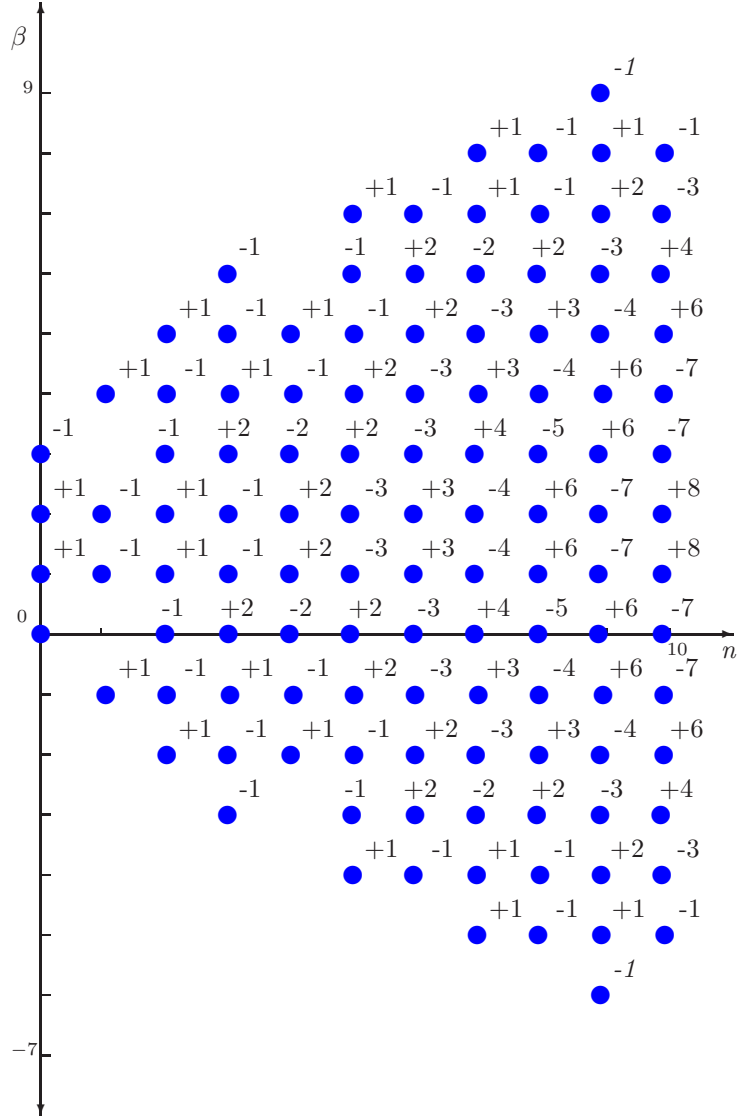
\begin{picture}(120.00,200.00)(-15.00,-09.00)
\put(00.00,80.00){\vector(0,50){105.00}}
\put(00.00,80.00){\vector(0,-50){80.00}}
\put(00.00,80.00){\vector(1,0){115.00}}
\multiput(10.00,80.00)(10.50,00.00){10}{\line(0,1){1.00}}
\put(104.00,76.50){\scriptsize $10$}
\put(113.00,76.00){n}
\put(-04.00,82.00){\scriptsize $0$}
\multiput(00.00,10.00)(00.00,10.00){17}{\line(1,0){1.50}}
\put(-05.00,10.00){\scriptsize $-7$}
\put(-05.00,178.00){$\beta$}
\put(-03.00,170.00){\scriptsize $9$}
\multiput(00.00,10.00)(00.00,10.00){17}{\line(1,0){1.50}}
\put(00.00,80.00){\usebox{\fstring}}
\put(02.00,83.00){\usebox{\fstringzm}}
\put(00.00,110.00){\usebox{\fstring}}
\put(02.00,113.00){\usebox{\fstringm}}
\put(31.00,140.00){\usebox{\fsstring}}
\put(33.00,143.00){\usebox{\fsstringm}}
\put(93.00,170.00){{\color{blue}\circle*{3}}}
\put(95.00,173.00){-1}
\put(31.00,50.00){\usebox{\fsstring}}
\put(33.00,53.00){\usebox{\fsstringm}}
\put(93.00,20.00){{\color{blue}\circle*{3}}}
\put(95.00,23.00){-1}
\put(00.00,90.00){\usebox{\fstrings}}
\put(02.00,93.00){\usebox{\fstringsm}}
\put(00.00,100.00){\usebox{\fstrings}}
\put(02.00,103.00){\usebox{\fstringsm}}
\put(10.80,120.00){\usebox{\fsstrings}}
\put(12.80,123.00){\usebox{\fsstringsm}}
\put(20.90,130.00){\usebox{\fssstrings}}
\put(22.90,133.00){\usebox{\fssstringsm}}
\put(10.80,70.00){\usebox{\fsstrings}}
\put(12.80,73.00){\usebox{\fsstringsm}}
\put(20.90,60.00){\usebox{\fssstrings}}
\put(22.90,63.00){\usebox{\fssstringsm}}
\put(51.80,150.00){\usebox{\fsssstrings}}
\put(53.80,153.00){\usebox{\fsssstringsm}}
\put(51.80,40.00){\usebox{\fsssstrings}}
\put(53.80,43.00){\usebox{\fsssstringsm}}
\put(72.50,160.00){\usebox{\fssssstrings}}
\put(74.50,163.00){\usebox{\fssssstringsm}}
\put(72.50,30.00){\usebox{\fssssstrings}}
\put(74.50,33.00){\usebox{\fssssstringsm}}
\end{picture}

Now we are able to construct the branching rules and explicitly reduce the
module $L^{\omega _{0}}$ with respect to the subalgebra $A_{2}^{(2)}$.
Remember that in terms of $A_{2}^{(2)}$ the diagram $\mathcal{N}^{\omega
_{0}}$ is located in the subspace of level $k=2$. Applying the formula (\ref
{recurrent relation}) in the sublattice with $n=0$ we find the first
nontrivial value for the weight with the highest (for $n=0$) anomalous
multiplicity , that is for the vector $\left( 1,2,0\right) $ where the
branching coefficient is evident, $k_{\left( 1,2,0\right) }^{\left( \mu
\right) }=1$. Implementing the recurrence procedure we obtain the set $%
\widehat{\Psi _{A_{2}^{(2)}}}$ of singular weights:
\begin{picture}(120.00,200.00)(-05.00,05.00)
\put(120.00,100.00){\vector(0,50){100.00}}
\put(120.00,100.00){\vector(0,-50){120.00}}
\put(120.00,100.00){\vector(-1,0){110.00}}
\put(125.00,110.00){\line(1,0){20.00}}
\put(125.00,120.00){\line(1,0){20.00}}
\put(135.00,113.00){$C_{A^{\left( 2\right)}_2}$}
\put(05.00,110.00){\line(1,0){10.00}}
\put(05.00,120.00){\line(1,0){20.00}}
\multiput(20.00,100.00)(10.00,00.00){10}{\line(0,1){2.00}}
\put(17.00,97.00){\scriptsize $-10$}
\put(10.00,95.00){ n}
\put(126.00,97.00){\scriptsize $0$}
\multiput(120.00,-10.00)(00.00,10.00){21}{\line(1,0){1.50}}
\put(124.00,190.00){\scriptsize $9$}
\put(123.00,200.00){$\beta$}
\put(124.00,00.00){\scriptsize $-10$}
\put(120.00,100.00){\circle*{2}}
\put(18.00,97.00){\usebox{\fbranch}}
\put(20.00,102.00){\usebox{\fbranchm}}
\put(18.00,107.00){\usebox{\fbranch}}
\put(20.00,111.50){\usebox{\fbranchp}}
\put(18.00,147.00){\usebox{\fbrranch}}
\put(20.00,152.00){\usebox{\fbrranchm}}
\put(18.00,57.00){\usebox{\fbrranch}}
\put(20.00,62.50){\usebox{\fbrranchp}}
\put(08.00,157.00){\usebox{\fbrranch}}
\put(10.00,162.50){\usebox{\fbrranchp}}
\put(08.00,47.00){\usebox{\fbrranch}}
\put(10.00,52.50){\usebox{\fbrranchm}}
\put(21.00,117.00){\usebox{\fbranchs}}
\put(24.00,122.00){\usebox{\fbranchps}}
\put(21.00,87.00){\usebox{\fbranchs}}
\put(25.00,92.00){\usebox{\fbranchms}}
\put(11.00,67.00){\usebox{\fbranchs}}
\put(14.00,72.00){\usebox{\fbranchps}}
\put(11.00,137.00){\usebox{\fbranchs}}
\put(14.00,142.00){\usebox{\fbranchms}}
\put(20.00,167.00){\usebox{\fbrranchs}}
\put(24.00,172.00){\usebox{\fbrranchps}}
\put(20.00,37.00){\usebox{\fbrranchs}}
\put(24.00,42.00){\usebox{\fbrranchms}}
\put(19.00,187.00){\usebox{\fbrrranchs}}
\put(22.00,192.00){\usebox{\fbrrranchms}}
\put(19.00,17.00){\usebox{\fbrrranchs}}
\put(22.00,22.00){\usebox{\fbrrranchps}}
\put(116.00,101.80){\scriptsize $-2$}
\put(120.00,110.00){{\color{blue}\circle*{2}}}
\put(116.00,111.80){\scriptsize $+1$}
\put(120.00,80.00){{\color{blue}\circle*{2}}}
\put(116.00,81.80){\scriptsize $+2$}
\put(120.00,70.00){{\color{blue}\circle*{2}}}
\put(116.00,71.80){\scriptsize $-1$}
\put(100.00,120.00){{\color{blue}\circle*{2}}}
\put(96.00,121.80){\scriptsize $+2$}
\put(100.00,150.00){{\color{blue}\circle*{2}}}
\put(96.00,151.80){\scriptsize $-1$}
\put(100.00,60.00){{\color{blue}\circle*{2}}}
\put(96.00,61.80){\scriptsize $-2$}
\put(100.00,30.00){{\color{blue}\circle*{2}}}
\put(96.00,31.80){\scriptsize $+1$}
\put(90.00,160.00){{\color{blue}\circle*{2}}}
\put(86.00,161.80){\scriptsize $+2$}
\put(90.00,140.00){{\color{blue}\circle*{2}}}
\put(86.00,141.80){\scriptsize $-2$}
\put(90.00,110.00){{\color{blue}\circle*{2}}}
\put(86.00,111.80){\scriptsize $-2$}
\put(90.00,70.00){{\color{blue}\circle*{2}}}
\put(86.00,71.80){\scriptsize $+2$}
\put(90.00,40.00){{\color{blue}\circle*{2}}}
\put(86.00,41.80){\scriptsize $+2$}
\put(90.00,20.00){{\color{blue}\circle*{2}}}
\put(86.00,21.80){\scriptsize $-2$}
\put(70.00,180.00){{\color{blue}\circle*{2}}}
\put(66.00,181.80){\scriptsize $-2$}
\put(70.00,150.00){{\color{blue}\circle*{2}}}
\put(66.00,151.80){\scriptsize $+2$}
\put(70.00,120.00){{\color{blue}\circle*{2}}}
\put(66.00,121.80){\scriptsize $+2$}
\put(70.00,60.00){{\color{blue}\circle*{2}}}
\put(66.00,61.80){\scriptsize $-2$}
\put(70.00,30.00){{\color{blue}\circle*{2}}}
\put(66.00,31.80){\scriptsize $-2$}
\put(70.00,00.00){{\color{blue}\circle*{2}}}
\put(66.00,01.80){\scriptsize $+2$}
\put(60.00,190.00){{\color{blue}\circle*{2}}}
\put(56.00,191.80){\scriptsize $+1$}
\put(60.00,140.00){{\color{blue}\circle*{2}}}
\put(56.00,141.80){\scriptsize $-2$}
\put(60.00,40.00){{\color{blue}\circle*{2}}}
\put(56.00,41.80){\scriptsize $+2$}
\put(60.00,-10.00){{\color{blue}\circle*{2}}}
\put(56.00,-08.50){\scriptsize $-1$}
\put(30.00,200.00){{\color{blue}\circle*{2}}}
\put(24.00,201.80){\scriptsize $+2$}
\put(30.00,190.00){{\color{blue}\circle*{2}}}
\put(26.00,191.80){\scriptsize $-2$}
\put(30.00,100.00){{\color{blue}\circle*{2}}}
\put(26.00,101.80){\scriptsize $-2$}
\put(30.00,80.00){{\color{blue}\circle*{2}}}
\put(26.00,81.80){\scriptsize $+2$}
\put(30.00,-10.00){{\color{blue}\circle*{2}}}
\put(26.00,-08.50){\scriptsize $+2$}
\put(29.00,197.00){{\Huge $\ast$}}
\put(32.00,202.00){-1}
\put(27.00,9.00){{\Huge $\ast$}}
\put(32.00,12.00){+1}
\end{picture}
\newpage We have performed the branching in terms of the singular weights $%
\widehat{\Psi _{A_{2}^{(2)}}^{\left( \xi \right) }}$ of the submodules $%
L_{A_{2}^{(2)}}^{\xi }$ in the decomposition $L_{A_{2}^{(1)}\downarrow
A_{2}^{(2)}}^{\mu }=\bigoplus\limits_{\xi \in P_{A_{2}^{(2)}}^{+}}b_{\xi
}^{\left( \mu \right) }L_{A_{2}^{(2)}}^{\xi }$ . Now it is quite easy to
extract the branching coefficients $b_{\xi }^{\left( \mu \right) }$. The
intersection
\begin{equation*}
\widehat{\Psi _{A_{2}^{(2)}}}\bigcap \overline{C}_{A_{2}^{(2)}}
\end{equation*}
gives the set of highest weights and their multiplicities $b_{\xi }^{\left(
\mu \right) }$ and the branching is
\begin{eqnarray*}
L_{A_{2}^{(1)}\downarrow A_{2}^{(2)}}^{\omega _{0}}
&=&L_{A_{2}^{(2)}}^{\omega _{0}}\left( 0\right) \oplus
L_{A_{2}^{(2)}}^{2\omega _{1}}\left( -1\right) \oplus
2L_{A_{2}^{(2)}}^{2\omega _{1}}\left( -3\right) \oplus
L_{A_{2}^{(2)}}^{\omega _{0}}\left( -4\right)  \\
&&\oplus 2L_{A_{2}^{(2)}}^{2\omega _{1}}\left( -5\right) \oplus
2L_{A_{2}^{(2)}}^{\omega _{0}}\left( -6\right) \oplus
4L_{A_{2}^{(2)}}^{2\omega _{1}}\left( -7\right)  \\
&&\oplus 3L_{A_{2}^{(2)}}^{\omega _{0}}\left( -8\right) \oplus
5L_{A_{2}^{(2)}}^{2\omega _{1}}\left( -9\right) \oplus
4L_{A_{2}^{(2)}}^{\omega _{0}}\left( -10\right) \oplus \ldots
\end{eqnarray*}
(Notice that as far as we have shifted the set $\Phi _{A_{2}^{(2)}\subset
A_{2}^{(1)}}$ the Weyl chamber $\overline{C}_{A_{2}^{(2)}}$ is also shifted
correspondingly.) The result can be presented in terms of two branching
functions
\begin{eqnarray*}
b_{I}^{\left( \mu \right) }\left( q\right)
&=&1+q^{4}+2q^{6}+3q^{8}+4q^{10}+\ldots  \\
b_{II}^{\left( \mu \right) }\left( q\right)
&=&q+2q^{3}+2q^{5}+4q^{7}+5q^{9}+\ldots
\end{eqnarray*}
\end{example}

\section{Conclusions}

We have demonstrated that the decompositions of integrable highest weight
modules of a simple Lie algebra (classical or affine) with respect to its
reductive subalgebra obey the (infinite) set of algebraic relations. These
relations originate from the properties of the singular vectors of the
module $L_{\frak{g}}$ considered as the highest weights of the Verma modules
$M_{\frak{a}}$. This gives
rise to the recursion relations for the branching coefficients.

The properties stated above are encoded in the subset $\Gamma _{\frak{g}
\supset \frak{a}}$ of the weight lattice $P_{\frak{a}}$ called the fan of
the injection. The fan depends only on the map $\frak{a} \longrightarrow
\frak{g}$. It describes the injection (whenever it is regular or special)
just as the root system describes the injection $\frak{h} \left( \frak{g}
\right) \longrightarrow \frak{g}$ of a Cartan subalgebra. Thus in the
simplest case, when $\frak{a} = \frak{h} \left( \frak{g} \right)$, the
recursion procedure produces the weight diagram of a module $L_{\frak{g}}$.

When applied to a reduction of highest weight modules the recursion
described by the fan provides a highly effective tool to obtain the explicit
values of branching coefficients.


\begin{thebibliography}{99}
\bibitem{DJKMO}  E Date, M Jimbo, A Kuniba, T Miwa and M Okado, ''One
dimensional configuration sums in vertex models and affine Lie algerba
characters'', Preprint RIMS-631, 1988; Lett. Math. Phys. 17, pp. 69-77, 1989

\bibitem{LD}  L D Faddeev, ''How Bethe ansatz works for integrable model'',
in: Quantum Symmetries/Symetries Quantique. Proc. Les Houches summer school,
LXIV, Eds. A Connes, K Gawedski, J Zinn-Justen; North-Holland, pp. 149-211,
1998; hep-th/9605187

\bibitem{KAZA}  V A Kazakov, K Zarembo, ''Classical/quantum integrability in
non-compact sector of AdS/CFT'', JHEP, 0410, 060, 2004

\bibitem{BE}  N Beisert, ''The Dilatation Operator of N=4 Super Y-M Theory
and Integrability'', Phys. Rep., 405, pp. 1-202, 2005; hep-th/0407277

\bibitem{Kac}  V Kac, ''Infinite Dimensional Lie Algebras'', (3rd edition),
Cambridge Univ. Press, 1990

\bibitem{Wak1}  M Wakimoto, ''Infinite-Dimensional Lie Algebras'',
Translations of Mathematical Monographs, Vol. 195, AMS, 2001

\bibitem{BGG}  I N Bernstein, I M Gelfand, S I Gelfand, ''Differential
operators on the basic affine space and a study of $\gamma $-modules, Lie
groups and their representations'' in ''Summer school of Bolyai Janos Math.
Soc.'', Ed. I M Gelfand, pp. 21-64, Halsted Press, NY, 1975

\bibitem{FauKing}  B Fauser, P D Jarvis, R C King and B G Wybourn, "New
branching rules induced by plethysm", arXiv:math-RT/0308043, 2003

\bibitem{Hwang}  S Hwang and H Rhedin, "General branching functions of
affine Lie algebras", arXiv:hep-th/9408087, 1994

\bibitem{FeigJimbo}  B Feigin, E Feigin, M Jimbo, T Miwa and E Mukhin
"Principal $\widehat{\frak{sl_3}}$ subspaces and quantum Toda Hamiltonians",
arXiv:0707.1635v2, 2007

\bibitem{LF}  V D Lyakhovsky and S Yu Melnikov, "Recursion relations and
branching rules for simple Lie algebras", J. Phys. A, 29 (1996) 1075-1087.

\bibitem{LDu}  V D Lyakhovsky, ''Recurrent properties of affine Lie algebra
representations'', Supersymmetry and Quantum Symmetry, Dubna, August 2007,
to be published.

\bibitem{Ful} W. Fulton, J. Harris, "Representation theory. A first course",
Graduate text in Math. v.129, Springer, NY, 1991

\bibitem{Burbaki} N. Burbaki "Elements de Mathematique. Groupes et Algebres
de Lie", Hermann, 1975

\end{thebibliography}
\end{document}